\documentclass[12pt,letterpaper,oneside,reqno]{amsart}
\usepackage{amsfonts}
\usepackage{amsmath}
\usepackage{amssymb}
\usepackage{amsthm}
\usepackage{float}
\usepackage{mathrsfs}
\usepackage[font=small,labelfont=bf]{caption}
\usepackage[left=1in,right=1in,bottom=1in,top=1in]{geometry}
\usepackage[pdfpagelabels,hyperindex,colorlinks=true,linkcolor=blue,urlcolor=magenta,citecolor=green]{hyperref}
\usepackage{graphicx}
\usepackage{array}
\usepackage{etoolbox}
\apptocmd{\sloppy}{\hbadness 10000\relax}{}{}
\newcommand \convPower [3]{(#1 \sp{#2} \ast #1 \sp{#2}) [#3]}
\newcommand \bernoulli [2][B] {{#1}\sb{#2}}
\newcommand \multifoldSum [2][x]{{#1}\sb{1} + {#1}\sb{2} + \cdots + {#1}\sb{#2}}

\newcommand \coeffA [3][A] {{\mathbf{#1}} \sb{#2,#3}}
\newcommand \coeffH [4][H] {{\mathbf{#1}} \sb{#2,#3} (#4)}
\newcommand \polynomialX [4][X] {{\mathbf{#1}}\sb{#2,#3} (#4)}
\newcommand \polynomialP [4][P]{{\mathbf{#1}}\sp{#2} \sb{#3}(#4)}
\newcommand \polynomialL [4][L]{{\mathbf{#1}}\sb{#2}(#3,#4)}

\let\svthefootnote\thefootnote
\newcommand\freefootnote[1]{%
    \let\thefootnote\relax%
    \footnotetext{#1}%
    \let\thefootnote\svthefootnote%
}

\newtheorem{thm}{Theorem}[section]
\newtheorem{cor}[thm]{Corollary}
\newtheorem{prop}[thm]{Proposition}
\newtheorem{lem}[thm]{Lemma}
\newtheorem{ppty}[thm]{Property}

\numberwithin{equation}{section}

\title[On the link between binomial theorem and discrete convolution]
{On the link between binomial theorem and discrete convolution}
\author[Petro Kolosov]{Petro Kolosov}
\email{kolosovp94@gmail.com}
\urladdr{https://kolosovpetro.github.io}
\date{\today}
\hypersetup{
    pdftitle={On the link between binomial theorem and discrete convolution},
    pdfsubject={
        Binomial theorem,
        Discrete convolution,
        Power function,
        Polynomials,
        Convolution,
        Multinomial theorem,
        Binomial coefficients,
        Bernoulli numbers,
        Pascal triangle,
        Faulhaber's formula,
        Power sums,
        Worpitzky identity
    },
    pdfauthor={Petro Kolosov},
    pdfkeywords={
        Binomial theorem,
        Discrete convolution,
        Power function,
        Polynomials,
        Convolution,
        Multinomial theorem,
        Binomial coefficients,
        Bernoulli number,
        Pascal's triangle,
        Faulhaber's formula,
        Power sums,
        Worpitzky identity
    }
}
\begin{document}
    \begin{abstract}
        Let $\mathbf{P}^{m}_{b}(x)$ be a $2m+1$-degree polynomial in $x$ and $b \in \mathbb{R}$
\[
    \mathbf{P}^{m}_{b}(x) = \sum_{k=0}^{b-1} \sum_{r=0}^{m} \mathbf{A}_{m,r} k^r (x-k)^r
\]
where $\mathbf{A}_{m,r}$ are real coefficients.
In this manuscript, we introduce the polynomial $\mathbf{P}^{m}_{b}(x)$ and study its properties,
establishing a polynomial identity for odd-powers in terms of this polynomial.
Based on mentioned polynomial identity for odd-powers,
we explore the connection between the Binomial theorem and discrete convolution of odd-powers,
further extending this relation to the multinomial case.
All findings are verified using Mathematica programs.

    \end{abstract}

    \maketitle

    \tableofcontents

    \freefootnote{Sources: \url{https://github.com/kolosovpetro/OnTheBinomialTheoremAndDiscreteConvolution}}

    \section{Definitions} \label{sec:definitions-notations-and-conventions}
    We now set the following notation, which remains fixed for the remainder of this manuscript
\begin{itemize}
    \setlength\itemsep{1.6em}
    \item $\coeffA{m}{r}$ is a real coefficient defined recursively
    \begin{equation}
        \label{eq:def_coeff_a}
        \coeffA{m}{r} =
        \begin{cases}
        (2r+1)
            \binom{2r}{r} & \text{if } r=m \\
            (2r+1) \binom{2r}{r} \sum_{d=2r+1}^{m} \coeffA{m}{d} \binom{d}{2r+1} \frac{(-1)^{d-1}}{d-r}
            \bernoulli{2d-2r} & \text{if } 0 \leq r<m \\
            0 & \text{if } r<0 \text{ or } r>m
        \end{cases}
    \end{equation}
    where $m$ is non-negative integer and $\bernoulli{t}$ are Bernoulli numbers~\cite{WeissteinBernoulli}.
    It is assumed that $\bernoulli{1}=\frac{1}{2}$.

    \item $\polynomialP{m}{b}{x}$ is a $2m+1$-degree polynomial in $b,x\in\mathbb{R}$
    \begin{equation}
        \label{eq:def_polynomial_p}
        \polynomialP{m}{b}{x} = \sum_{k=0}^{b-1} \sum_{r=0}^{m} \coeffA{m}{r} k^r(x-k)^r
    \end{equation}

    \item $\coeffH{m}{t}{b}$ is a polynomial defined as
    \begin{equation}
        \label{eq:def_coeff_h}
        \coeffH{m}{t}{b}
        = \sum_{j=t}^{m} \binom{j}{t} \coeffA{m}{j} \frac{(-1)^j}{2j-t+1} \binom{2j-t+1}{b} \bernoulli{2j-t+1-b}
    \end{equation}
    integers $m,t,b$.

    \item $\polynomialX{m}{t}{j}$ is polynomial of degree $2m+1-t$ in $j\in\mathbb{R}$
    \begin{equation}
        \label{eq:def_coeff_x}
        \polynomialX{m}{t}{j} = (-1)^m \sum_{k=1}^{2m+1-t} \coeffH{m}{t}{k} \cdot j^k
    \end{equation}
    integers $m,t$.

    \item $\polynomialL{m}{x}{k}$ is $2m$ degree polynomial in $x,k\in\mathbb{R}$
    \begin{equation}
        \label{eq:def_polynomial_l}
        \polynomialL{m}{x}{k} = \sum_{r=0}^{m} \coeffA{m}{r} k^r(x-k)^r
    \end{equation}

    \item $(f\ast f)[n]$ is discrete convolution~\cite{damelin_discrete_convolution} of function $f$ defined over set of integers $\mathbb{Z}$
    \begin{align*}
    (f\ast f)[n]
        = \sum_{k} f(k) f(n-k)
    \end{align*}
    and its partial case for polynomials $n^j, \; n\geq a \in \mathbb{R}$
    \begin{align*}
        \convPower{n}{j}{x} =\sum_{k} k^j (x-k)^j [k\geq a][x-k\geq a] =\sum_{k=a}^{x-a} k^j (x-k)^j
    \end{align*}
\end{itemize}

    \clearpage

    \section{Introduction and main results} \label{sec:introduction}
    The polynomial $\polynomialP{m}{b}{x}$ is a $2m+1$-degree polynomial in $x,b\in\mathbb{R}$ defined as
\begin{align*}
    \polynomialP{m}{b}{x} = \sum_{k=0}^{b-1} \sum_{r=0}^{m} \coeffA{m}{r} k^r(x-k)^r
\end{align*}
where $\coeffA{m}{r}$ is a real coefficient.
By means of Lemma~\eqref{lemma_polynomial_p_and_odd_power},
the polynomial $\polynomialP{m}{b}{x}$ has the following relation with Binomial theorem~\cite{AbraSteg72}
\begin{align*}
    \polynomialP{m}{x+y}{x+y} = \sum_{r=0}^{2m+1} \binom{2m+1}{r} x^{2m+1-r} y^r
\end{align*}
On the other hand, polynomial $\polynomialP{m}{b}{x}$ might be expressed in terms of discrete convolution
of polynomial $n^j$.
For every $n\geq 0$
\begin{align*}
    \polynomialP{m}{x+1}{x} = \sum_{r=0}^{m} \coeffA{m}{r} \convPower{n}{r}{x}
\end{align*}
It is important to notice that  $n^r$ of discrete convolution $\convPower{n}{r}{x}$ evaluated at $x$
is implicit piecewise-defined polynomial such as
\begin{equation*}
    n^{r} =
    \begin{cases}
        \underbrace{n \cdot n \cdots n}_{\mathrm{r \; times}}, & \mbox{if } n \geq 0 \\
        0, & \mbox{otherwise}
    \end{cases}
\end{equation*}
Therefore, it is easy to notice the following identities in terms of Binomial theorem and discrete convolution,
see the corollaries~\eqref{cor_bin_exp_and_macaulay_conv} and~\eqref{cor_bin_exp_and_macaulay_conv_strict}.
For every $n \geq 0$
\begin{equation*}
    \sum_{r=0}^{m} \coeffA{m}{r} \convPower{n}{r}{x+y}
    = 1 + \sum_{r=0}^{2m+1} \binom{2m+1}{r} x^{2m+1-r} y^r
\end{equation*}
For every $n > 0$
\begin{equation*}
    \sum_{r=0}^{m} \coeffA{m}{r} \convPower{n}{r}{x+y}
    = -1 + \sum_{r=0}^{2m+1} \binom{2m+1}{r} x^{2m+1-r} y^r
\end{equation*}
Additionally, the following generalizations for the multinomial case are discussed in
the corollaries~\eqref{cor_mult_exp_and_macaulay_conv} and ~\eqref{cor_mult_exp_and_macaulay_conv_strict}.
For every $n \geq 0$
\begin{align*}
    \sum_{r=0}^{m} \coeffA{m}{r} \convPower{n}{r}{\multifoldSum{t}} =
    1 + \sum_{\multifoldSum[k]{t}=2m+1} \binom{2m+1}{k_1, k_2,\ldots, k_t} \prod_{\ell=1}^{t} x_\ell^{k_\ell}
\end{align*}
For every $n>0$
\begin{align*}
    \sum_{r=0}^{m} \coeffA{m}{r} \convPower{n}{r}{\multifoldSum{t}} =
    -1 + \sum_{\multifoldSum[k]{t}=2m+1} \binom{2m+1}{k_1, k_2,\ldots, k_t} \prod_{\ell=1}^{t} x_\ell^{k_\ell}
\end{align*}
A few polynomial identities are straightforward by means of
the theorems~\eqref{thm_odd_power_by_macaulays_convolution},~\eqref{thm_odd_power_by_macaulays_convolution_strict}.
More precisely, by the theorem~\eqref{thm_odd_power_by_macaulays_convolution} we have an odd-power identity as follows
\begin{equation*}
    x^{2m+1} = \sum_{r=0}^{m} \coeffA{m}{r} \sum_{k=0}^{x-1} k^r (x-k)^r
\end{equation*}
so that
\begin{align*}
    1 + x^{2m+1} = \sum_{r=0}^{m} \coeffA{m}{r} \convPower{n}{r}{x}
    = \sum_{r=0}^{m} \coeffA{m}{r} \sum_{k=0}^{x} k^r (x-k)^r
\end{align*}
From the other side, the theorem~\eqref{thm_odd_power_by_macaulays_convolution_strict} provides an odd-power
polynomial identity as follows
\begin{equation*}
    -1 + x^{2m+1} = \sum_{r=0}^{m} \coeffA{m}{r} \convPower{n}{r}{x}
    = \sum_{r=0}^{m} \coeffA{m}{r} \sum_{k=1}^{x-1} k^r (x-k)^r
\end{equation*}
For example,
\begin{align*}
    x^3 &= \sum_{k=1}^{x} 6k (x-k) + 1 \\
    x^5 &= \sum_{k=1}^{x} 30k^2 (x-k)^2 + 1 \\
    x^7 &= \sum_{k=1}^{x} 140 k^3 (x-k)^3 - 14k(x-k) + 1 \\
    x^9 &= \sum_{k=1}^{x} 630 k^4(x-k)^4 - 120k(x-k) + 1 \\
    x^{11} &= \sum_{k=1}^{x} 2772 k^5 (x-k)^5 + 660 k^2(x-k)^2 - 1386k(x-k) + 1 \\
    x^{13} &= \sum_{k=1}^{x} 51480 k^7 (x-k)^7 - 60060 k^3 (x-k)^3 + 491400 k^2 (x-k)^{2} - 450054 k (x-k) + 1 \\
\end{align*}
Moreover, the following binomials in terms of discrete convolution of polynomial $n^j$ are found,
see the equations~\eqref{eq:parametric-identity} and~\eqref{eq:parametric-identity-strict}.
For every $n \geq 0$
\begin{equation*}
    \begin{split}
    (x-2a)
        ^{2m+1} + 1 &= \sum_{r=0}^{m} \coeffA{m}{r} ((t-k)^r \ast (t-k)^r)[x] \\
                    &= \sum_{r=0}^{m} \coeffA{m}{r} \sum_{k=a}^{x-a} (k-a)^r (x-k-a)^r
    \end{split}
\end{equation*}
Similarly, the following binomial holds.
For every $n > 0$
\begin{equation*}
    \begin{split}
    (x-2a)
        ^{2m+1} - 1 &= \sum_{r=0}^{m} \coeffA{m}{r} ((t-k)^r \ast (t-k)^r)[x] \\
                    &= \sum_{r=0}^{m} \coeffA{m}{r} \sum_{k=a+1}^{x-a-1} (k-a)^r (x-k-a)^r
    \end{split}
\end{equation*}
This manuscript does not contain any historical context about the polynomial $\polynomialP{m}{b}{x}$,
either how exactly it was derived including all the milestones.
To get more information about the history of polynomials $\polynomialP{m}{b}{x}$, the reader can refer
to the manuscript
\begin{center}
    \href{https://kolosovpetro.github.io/pdf/HistoryAndOverviewOfPolynomialP.pdf}
    {\texttt{kolosovpetro.github.io/pdf/HistoryAndOverviewOfPolynomialP.pdf}}
\end{center}

    \section{Polynomial \texorpdfstring{$\polynomialP{m}{b}{x}$}{P[m,b,x]} and its properties}
    \label{sec:polynomial-p-and-its-properties}
    \label{sec:polynomial-p-and-their-properties}
We continue our mathematical journey from the short overview
of polynomial $\polynomialL{m}{x}{k}$ which is
an essential part of polynomial $\polynomialP{m}{b}{x}$ since that
$\polynomialP{m}{b}{x} = \sum_{k=0}^{b-1} \polynomialL{m}{x}{k}$.
Polynomial $\polynomialL{m}{x}{k}$ is a polynomial of degree $2m$ in $x,k\in\mathbb{R}$,
see definition~\eqref{eq:def_polynomial_l}.
In its explicit form the polynomial $\polynomialL{m}{x}{k}$ is as follows
\begin{equation*}
    \polynomialL{m}{x}{k} =
    \coeffA{m}{m} k^m(x-k)^m +
    \coeffA{m}{m-1} k^{m-1}(x-k)^{m-1} +
    \cdots +
    \coeffA{m}{0}
\end{equation*}
where $\coeffA{m}{r}$ are real coefficients defined by~\eqref{eq:def_coeff_a}.
Coefficients $\coeffA{m}{r}$ are nonzero for $r$ only within the range $r \in \{m\} \cup \left[0,\frac{m-1}{2}\right]$.
For example,
\begin{table}[H]
    \setlength\extrarowheight{-6pt}
    \begin{tabular}{c|cccccccc}
        $m/r$ & 0 & 1       & 2      & 3      & 4   & 5    & 6     & 7     \\
        \hline
        0     & 1 &         &        &        &     &      &       &       \\
        1     & 1 & 6       &        &        &     &      &       &       \\
        2     & 1 & 0       & 30     &        &     &      &       &       \\
        3     & 1 & -14     & 0      & 140    &     &      &       &       \\
        4     & 1 & -120    & 0      & 0      & 630 &      &       &       \\
        5     & 1 & -1386   & 660    & 0      & 0   & 2772 &       &       \\
        6     & 1 & -21840  & 18018  & 0      & 0   & 0    & 12012 &       \\
        7     & 1 & -450054 & 491400 & -60060 & 0   & 0    & 0     & 51480
    \end{tabular}
    \caption{Coefficients $\coeffA{m}{r}$. See the OEIS entries
    \href{https://oeis.org/A302971}{\texttt{A302971}},
        \href{https://oeis.org/A304042}{\texttt{A304042}}: \cite{kolosov2018numerator, kolosov2018denominator}.}
    \label{tab:table_of_coefficients_a}
\end{table}
Thus, the polynomial $\polynomialL{m}{x}{k}$ may also be written as
\begin{equation*}
    \polynomialL{m}{x}{k} = \coeffA{m}{m} k^m (x-k)^m + \sum_{r=0}^{\frac{m-1}{2}} \coeffA{m}{r} k^r (x-k)^r
\end{equation*}
For example, the polynomials $\polynomialL{m}{x}{k}$ for $0\leq m\leq 3$ are
\begin{equation*}
    \begin{split}
        \polynomialL{0}{x}{k}
        &= 1 \\
        \polynomialL{1}{x}{k}
        &= 6 k (x-k) + 1
        = -6 k^2 + 6 k x + 1 \\
        \polynomialL{2}{x}{k}
        &=30 k^2 (x-k)^2+1
        =30 k^4-60 k^3 x+30 k^2 x^2+1 \\
        \polynomialL{3}{x}{k}
        &= 140 k^3 (x-k)^3-14 k (x-k)+1 \\
        &=-140 k^6+420 k^5 x-420 k^4 x^2+140 k^3 x^3+14 k^2-14 k x+1
    \end{split}
\end{equation*}
It is important to notice that $\polynomialL{m}{x}{k}$ is symmetric over $x$
\begin{ppty}
    \label{ppty_symmetry_of_polynomial_l}
    For every $x,k\in\mathbb{R}$
    \begin{equation*}
        \polynomialL{m}{x}{k} = \polynomialL{m}{x}{x-k}
    \end{equation*}
\end{ppty}
This might be seen from the following tables
\begin{table}[H]
    \setlength\extrarowheight{-6pt}
    \begin{tabular}{c|cccccccc}
        $x/k$ & 0 & 1  & 2  & 3  & 4  & 5  & 6  & 7 \\
        \hline
        0     & 1 &    &    &    &    &    &    &   \\
        1     & 1 & 1  &    &    &    &    &    &   \\
        2     & 1 & 7  & 1  &    &    &    &    &   \\
        3     & 1 & 13 & 13 & 1  &    &    &    &   \\
        4     & 1 & 19 & 25 & 19 & 1  &    &    &   \\
        5     & 1 & 25 & 37 & 37 & 25 & 1  &    &   \\
        6     & 1 & 31 & 49 & 55 & 49 & 31 & 1  &   \\
        7     & 1 & 37 & 61 & 73 & 73 & 61 & 37 & 1
    \end{tabular}
    ~\caption{Values of $\polynomialL{1}{x}{k}$.
    See the OEIS entry \href{https://oeis.org/A287326}{\texttt{A287326}}, \cite{kolosov2017third}.}
    \label{tab:fig_1}
\end{table}
Another case, given $m=2$ we have the following values of $\polynomialL{2}{x}{k}$
\begin{table}[H]
    \setlength\extrarowheight{-6pt}
    \begin{tabular}{c|cccccccc}
        $x/k$ & 0 & 1    & 2    & 3    & 4    & 5    & 6    & 7 \\
        \hline
        0     & 1 &      &      &      &      &      &      &   \\
        1     & 1 & 1    &      &      &      &      &      &   \\
        2     & 1 & 31   & 1    &      &      &      &      &   \\
        3     & 1 & 121  & 121  & 1    &      &      &      &   \\
        4     & 1 & 271  & 481  & 271  & 1    &      &      &   \\
        5     & 1 & 481  & 1081 & 1081 & 481  & 1    &      &   \\
        6     & 1 & 751  & 1921 & 2431 & 1921 & 751  & 1    &   \\
        7     & 1 & 1081 & 3001 & 4321 & 4321 & 3001 & 1081 & 1
    \end{tabular}
    \caption{Values of $\polynomialL{2}{x}{k}$.
    See the OEIS entry \href{https://oeis.org/A300656}{\texttt{A300656}}, ~\cite{kolosov2018fifth}.}
    \label{tab:row-sums-give-fifth-power}
\end{table}
Note that row sums of the table~\eqref{tab:fig_1} are cubes of $x$.
Next we discuss the polynomial $\polynomialP{m}{b}{x}$.
In its extended form, the polynomial $\polynomialP{m}{b}{x}$ is
\begin{equation*}
    \begin{split}
        \polynomialP{m}{b}{x} = \sum_{k=0}^{b-1} \polynomialL{m}{x}{k}
        =\sum_{k=0}^{b-1} \sum_{r=0}^{m} \coeffA{m}{r} k^r(x-k)^r
        =\sum_{r=0}^{m} \coeffA{m}{r} \sum_{k=0}^{b-1} k^r(x-k)^r
    \end{split}
\end{equation*}
By means of binomial theorem $(x-y)^n = \sum_{k=0}^{n} (-1)^{k} \binom{n}{k} x^{n-k} y^{k}$,
\begin{equation*}
    \begin{split}
        \polynomialP{m}{b}{x}
        &=\sum_{r=0}^{m} \coeffA{m}{r} \sum_{k=0}^{b-1} k^r \sum_{j=0}^{r} (-1)^{j} \binom{r}{j} x^{r-j} k^{j} \\
        &=\sum_{r=0}^{m} \coeffA{m}{r} \sum_{k=0}^{b-1} \sum_{j=0}^{r} (-1)^{j} \binom{r}{j} x^{r-j} k^{r+j} \\
        &=\sum_{r=0}^{m} \coeffA{m}{r} \sum_{j=0}^{r} (-1)^{j} x^{r-j} \binom{r}{j} \sum_{k=0}^{b-1} k^{r+j} \\
    \end{split}
\end{equation*}
However, by the symmetry~\eqref{ppty_symmetry_of_polynomial_l} of $\polynomialL{m}{x}{k}$ the polynomial
$\polynomialP{m}{b}{x}$ may also be written in the form
\begin{equation*}
    \begin{split}
        \polynomialP{m}{b}{x}
        &=\sum_{k=1}^{b} \sum_{r=0}^{m} \coeffA{m}{r} k^r(x-k)^r
        =\sum_{k=1}^{b} \sum_{r=0}^{m} \coeffA{m}{r} k^r \sum_{t=0}^{r} (-1)^{r-t} x^t \binom{r}{t} k^{r-t} \\
        &=\sum_{t=0}^{m} x^t
        \underbrace{\sum_{k=1}^{b} \sum_{r=t}^{m} (-1)^{r-t} \binom{r}{t} \coeffA{m}{r} k^{2r-t}}_{(-1)^{m-t} \polynomialX{m}{t}{b}}
    \end{split}
\end{equation*}
Note that
$\sum_{k=1}^{b} \sum_{r=t}^{m} (-1)^{r-t} \binom{r}{t} \coeffA{m}{r} k^{2r-t}$
is the
$(-1)^{m-t} \polynomialX{m}{t}{b}$.
From this formula it may be not immediately clear why $\polynomialX{m}{t}{b}$ represent polynomials in $b$.
However, this can be seen if we change the summation order and use Faulhaber's formula
$\sum_{k=1}^{n} k^{p}=\frac{1}{p+1}\sum _{j=0}^{p} \binom{p+1}{j} \bernoulli{j} n^{p+1-j}$
to obtain
\begin{equation*}
    \polynomialX{m}{t}{b} = (-1)^m \sum_{r=t}^{m} \binom{r}{t} \coeffA{m}{r} \frac{(-1)^r}{2r-t+1}
    \sum_{\ell=0}^{2r-t} \binom{2r-t+1}{\ell} \bernoulli{\ell} b^{2r-t+1-\ell}
\end{equation*}
Introducing $k=2r-t+1-\ell$ we further get the formula
\begin{equation*}
    \polynomialX{m}{t}{b} = (-1)^m \sum_{k=1}^{2m-t+1} b^k
    \underbrace{\sum_{r=t}^m \binom{r}{t} \coeffA{m}{r} \frac{(-1)^r}{2r-t+1} \binom{2r-t+1}{k}
    \bernoulli{2r-t+1-k}}_{\coeffH{m}{t}{k}}
\end{equation*}
Polynomials $\polynomialX{3}{t}{b}, \; 0\leq t \leq 3$ are
\begin{equation*}
    \begin{split}
        \polynomialX{3}{0}{j}
        &=7 b^2 - 28 b^3 + 70 b^5 - 70 b^6 + 20 b^7 \\
        \polynomialX{3}{1}{j}
        &=7 b - 42 b^2 + 175 b^4 - 210 b^5 + 70 b^6 \\
        \polynomialX{3}{2}{j}
        &=-14 b + 140 b^3 - 210 b^4 + 84 b^5 \\
        \polynomialX{3}{3}{j}
        &=35 b^2 - 70 b^3 + 35 b^4
    \end{split}
\end{equation*}
Polynomials $\coeffH{3}{t}{k}$ are defined by~\eqref{eq:def_coeff_h} and examples for $m=3, \; 0\leq t \leq 3$ are
\begin{equation*}
    \begin{split}
        \coeffH{3}{0}{k}
        &=\bernoulli{1-k} \binom{1}{k} + \frac{14}{3} \bernoulli{3-k} \binom{3}{k} - 20 \bernoulli{7 - k} \binom{7}{k} \\
        \coeffH{3}{1}{k}
        &=7 \bernoulli{2-k} \binom{2}{k} - 70 \bernoulli{6-k} \binom{6}{k} \\
        \coeffH{3}{2}{k}
        &=-84 \bernoulli{5-k} \binom{5}{k} \\
        \coeffH{3}{3}{k}
        &=-35 \bernoulli{4-k} \binom{4}{k}
    \end{split}
\end{equation*}
It gives us an opportunity to overview the polynomial $\polynomialP{m}{b}{x}$ from the different prospective,
for instance
\begin{equation}
    \label{eq:p_all_forms}
    \polynomialP{m}{b}{x}
    =\sum_{r=0}^{m} (-1)^{m-r} \polynomialX{m}{r}{b} \cdot x^r
    =\sum_{r=0}^{m} \sum_{\ell=1}^{2m-r+1} (-1)^{2m-r} \coeffH{m}{r}{\ell} \cdot b^\ell \cdot x^r
\end{equation}
Equation~\eqref{eq:p_all_forms} clearly states why $\polynomialP{m}{b}{x}$ is polynomial in $x,b$.
For example,
\begin{align*}
    \polynomialP{0}{b}{x}
    &=b \\
    \polynomialP{1}{b}{x}
    &=3 b^2 - 2 b^3 - 3 b x + 3 b^2 x \\
    \polynomialP{2}{b}{x}
    &=10 b^3 - 15 b^4 + 6 b^5 - 15 b^2 x + 30 b^3 x - 15 b^4 x + 5 b x^2 - 15 b^2 x^2 + 10 b^3 x^2 \\
    \polynomialP{3}{b}{x}
    &=-7 b^2 + 28 b^3 - 70 b^5 + 70 b^6 - 20 b^7
    + 7 b x - 42 b^2 x + 175 b^4 x - 210 b^5 x + 70 b^6 x \\
    &+ 14 b x^2 - 140 b^3 x^2 + 210 b^4 x^2 - 84 b^5 x^2
    + 35 b^2 x^3 - 70 b^3 x^3 + 35 b^4 x^3
\end{align*}
The following property is also true in terms of the polynomial $\polynomialP{m}{b}{x}$
\begin{ppty}
    \label{prop_p_identity}
    For every $m\in \mathbb{N}, \; x,b\in\mathbb{R}$
    \begin{equation*}
        \polynomialP{m}{b+1}{x} = \polynomialP{m}{b}{x} + \polynomialL{m}{x}{b}
    \end{equation*}
\end{ppty}

    \section{Relation between the polynomial \texorpdfstring{$\polynomialP{m}{b}{x}$}{P[m,b,x]} and Binomial theorem}
    \label{sec:odd-binomial-expansion-as-partial-case-of-polynomial-p}
    \begin{lem}
    \label{lemma_polynomial_p_and_odd_power}
    For every $m\in\mathbb{N}, \; x,y\in\mathbb{R}$
    \begin{equation*}
        \polynomialP{m}{x+y}{x+y} = \sum_{r=0}^{2m+1} \binom{2m+1}{r} x^{2m+1-r} y^r
    \end{equation*}
\end{lem}
By means of lemma~\ref{lemma_polynomial_p_and_odd_power} and equation~\eqref{eq:p_all_forms} the following
polynomial identities straightforward
\begin{equation*}
    x^{2m+1}
    =\sum_{r=0}^{m} \sum_{\ell=1}^{2m-r+1} (-1)^{2m-r} \coeffH{m}{r}{\ell} \cdot x^{\ell+r}
    =\sum_{r=0}^{m} (-1)^{m-r} \polynomialX{m}{r}{x} \cdot x^r
\end{equation*}
For instance,
\begin{equation*}
    \polynomialP{2}{x+y}{x+y} = (x + y) (x^4 + 4 x^3 y + 6 x^2 y^2 + 4 x y^3 + y^4).
\end{equation*}
In addition, the following identities hold
\begin{equation*}
    \begin{split}
    (x+y)
        ^{2m+1}
        &=\sum_{r=0}^{m} \sum_{\ell=1}^{2m-r+1} (-1)^{2m-r} \coeffH{m}{r}{\ell} \cdot (x+y)^{\ell+r} \\
        &=\sum_{r=0}^{m} (-1)^{m-r} \polynomialX{m}{r}{x+y} \cdot (x+y)^r
    \end{split}
\end{equation*}
Obviously, Multinomial expansion of $t$-fold sum $(\multifoldSum{t})^{2m+1}$ can be reached
by $\polynomialP{m}{b}{\multifoldSum{t}}$ as well
\begin{cor}
    For all $x_1,x_2,\ldots, x_t\in\mathbb{R}, \; m \in \mathbb{N}$
    \begin{equation*}
        \polynomialP{m}{\multifoldSum{t}}{\multifoldSum{t}}
        =
        \sum_{\multifoldSum[k]{t}=2m+1} \binom{2m+1}{k_1, k_2,\ldots, k_t} \prod_{s=1}^{t} x_t^{k_s}
    \end{equation*}
\end{cor}
Moreover, the following multinomial identities hold
\begin{equation*}
    \begin{split}
    (\multifoldSum{t})
        ^{2m+1}
        &=\sum_{r=0}^{m} \sum_{\ell=1}^{2m-r+1} (-1)^{2m-r} \coeffH{m}{r}{\ell} \cdot (\multifoldSum{t})^{\ell+r} \\
        &=\sum_{r=0}^{m} (-1)^{m-r} \polynomialX{m}{r}{\multifoldSum{t}} \cdot (\multifoldSum{t})^r
    \end{split}
\end{equation*}

    \section{Polynomial \texorpdfstring{$\polynomialP{m}{b}{x}$}{P[m,b,x]} in terms of Discrete convolution}
    \label{sec:relation-between-p-and-convolution-of-polynomials}
    In this section we discuss the relation between $\polynomialP{m}{b}{x}$ and discrete convolution of
polynomials.
To show that $\polynomialP{m}{b}{x}$ involves the discrete convolution of polynomial $n^r$
recall the definition of the polynomial $\polynomialP{m}{b}{x}$
\begin{equation*}
    \polynomialP{m}{b}{x} = \sum_{k=0}^{b-1} \sum_{r=0}^{m} \coeffA{m}{r} k^r (x-k)^r
    = \sum_{r=0}^{m} \coeffA{m}{r} \sum_{k=0}^{b-1} k^r (x-k)^r
\end{equation*}
A discrete convolution of defined over set of integers $\mathbb{Z}$ function $f$ is
\begin{equation*}
(f \ast f)[n]
    = \sum_{k} f(k) f(n-k)
\end{equation*}
General formula of discrete convolution for the polynomial $n^j, \; n\geq a \in \mathbb{R}$
can be derived immediately
\begin{align*}
    \convPower{n}{j}{x}
    &=\sum_{k} k^j (x-k)^j [k\geq a][x-k\geq a] \\
    &=\sum_{k} k^j (x-k)^j [k\geq a][k\leq x-a] \\
    &=\sum_{k} k^j (x-k)^j [a \leq k \leq x-a] \\
    &=\sum_{k=a}^{x-a} k^j (x-k)^j
\end{align*}
where $[a \leq k \leq x-a]$ is Iverson's bracket~\cite{iverson_apl, knuth_two_notes_on_notation}.
\begin{lem}
    \label{lemma_disc_conv_identity}
    For every $n\in\mathbb{N}, \; x\in\mathbb{R}$ and $n\geq 0$
    \[
        \convPower{n}{r}{x} = \sum_{k=0}^{x} k^r (x-k)^r
    \]
\end{lem}
It is of first importance to keep in mind that  $n^r$ of discrete convolution $\convPower{n}{r}{x}$ evaluated at $x$
is an implicit piecewise-defined polynomial such as
\begin{equation*}
    n^{r} =
    \begin{cases}
        \underbrace{n \cdot n \cdots n}_{\mathrm{r \; times}}, & \mbox{if } n \geq 0 \\
        0, & \mbox{otherwise}
    \end{cases}
\end{equation*}
Thus, the corollary follows
\begin{cor}
    \label{cor_polynomial_p_and_macaulay_convolution}
    By Lemma~\ref{lemma_disc_conv_identity} the polynomial $\polynomialP{m}{b}{n}$ might be expressed in terms
    of discrete convolution as follows, for every $n\geq 0$
    \begin{align*}
        \polynomialP{m}{x+1}{x} = \sum_{r=0}^{m} \coeffA{m}{r} \convPower{n}{r}{x}
    \end{align*}
\end{cor}
Therefore, another polynomial identity follows
\begin{thm}
    \label{thm_odd_power_by_macaulays_convolution}
    By Lemma~\ref{lemma_polynomial_p_and_odd_power}, Corollary~\ref{cor_polynomial_p_and_macaulay_convolution}
    and property~\ref{prop_p_identity}, for every $m\in\mathbb{N}, \; x\in\mathbb{R}$ and $n\geq 0$
    \begin{align*}
        1 + x^{2m+1} = \sum_{r=0}^{m} \coeffA{m}{r} \convPower{n}{r}{x}
        = \sum_{r=0}^{m} \coeffA{m}{r} \sum_{k=0}^{x} k^r (x-k)^r
    \end{align*}
\end{thm}
Now we notice the following identity in terms of polynomial $\polynomialP{m}{b}{x}$ and
discrete convolution $\convPower{n}{j}{x}$
\begin{prop}
    \label{prop_polynomial_p_and_macaulay_convolution_strict}
    For every $m \in \mathbb{N}, \; x\in\mathbb{R}$ and $n \geq 1$
    \begin{equation*}
        \begin{split}
            \polynomialP{m}{x}{x}
            &=\sum_{r=0}^{m} \coeffA{m}{r} \left(0^r x^r + \sum_{k=1}^{x-1} k^r (x-k)^r \right) \\
            &=\sum_{r=0}^{m} \coeffA{m}{r} 0^r x^r + \sum_{r=0}^{m} \coeffA{m}{r} \convPower{n}{r}{x} \\
            &=1 + \sum_{r=0}^{m} \coeffA{m}{r} \convPower{n}{r}{x}
        \end{split}
    \end{equation*}
\end{prop}
Since that for all $r$ in $\coeffA{m}{r} 0^r x^r$ we have
\begin{equation*}
    \coeffA{m}{r} 0^r x^r =
    \begin{cases}
        1, & \mbox{if } r=0 \\
        0, & \mbox{if } r>0
    \end{cases}
\end{equation*}
Above is true because $\coeffA{m}{0}=1$ for every $m\in\mathbb{N}$, and $x^0 = 1$
for every $x$, see~\cite{graham1994concrete}.
Hence, the following identity between $\polynomialP{m}{b}{x}$ and
discrete convolution $\convPower{n}{j}{x}$ holds
\begin{thm}
    \label{thm_odd_power_by_macaulays_convolution_strict}
    By Lemma~\ref{lemma_polynomial_p_and_odd_power} and
    Proposition~\ref{prop_polynomial_p_and_macaulay_convolution_strict},
    for every $m\in\mathbb{N}, \; x\in\mathbb{R}$ and $n > 0$
    \begin{equation*}
        -1 + x^{2m+1} = \sum_{r=0}^{m} \coeffA{m}{r} \convPower{n}{r}{x}
        = \sum_{r=0}^{m} \coeffA{m}{r} \sum_{k=1}^{x-1} k^r (x-k)^r
    \end{equation*}
\end{thm}
\begin{cor}
    \label{cor_sum_of_coeffs_a}
    By Theorem~\ref{thm_odd_power_by_macaulays_convolution_strict}, for all $m\in\mathbb{N}$
    \begin{equation*}
        \sum_{r=0}^{m} \coeffA{m}{r} = 2^{2m+1} - 1
    \end{equation*}
\end{cor}
Corollary~\ref{cor_sum_of_coeffs_a} holds since that convolution $\convPower{n}{j}{x}=1, \; n > 0$
for each $r$ and $x=2$.

    \section{Relation between Binomial theorem and Discrete convolution}
    \label{sec:relation-between-binomial-theorem-and-discrete-convolution}
    \begin{cor}
    \label{cor_bin_exp_and_macaulay_conv}
    (Generalization of Theorem~\ref{thm_odd_power_by_macaulays_convolution} for Binomials.)
    For every $m\in\mathbb{N}, \; x,y\in\mathbb{R}$ and $n\geq 0$
    \begin{equation*}
        \sum_{r=0}^{m} \coeffA{m}{r} \convPower{n}{r}{x+y}
        =
        1 + \sum_{r=0}^{2m+1} \binom{2m+1}{r} x^{2m+1-r} y^r
    \end{equation*}
\end{cor}
For example, given $m=0,1,2$ the Corollary~\ref{cor_bin_exp_and_macaulay_conv} yields
\begin{align*}
    \sum_{r=0}^{0} \coeffA{0}{r} \convPower{n}{r}{x+y}
    &= 1 + x + y \\
    \sum_{r=0}^{1} \coeffA{1}{r} \convPower{n}{r}{x+y}
    &= 1 + x + y - (x + y) (1 + x + y) (1 - 3 x - 3 y + 2 (x + y)) \\
    &= 1 + x^3 + 3 x^2 y + 3 x y^2 + y^3\\
    \sum_{r=0}^{2} \coeffA{2}{r} \convPower{n}{r}{x+y}
    &=1 + x + y + (x + y) (1 + x + y) \left(-1 + x + 5 x^2 + y + 10 x y + 5 y^2\right. \\
    &-15 x (x + y) + 10 x^2 (x + y) - 15 y (x + y) + 20 x y (x + y) \\
    &+ 10 y^2 (x + y) +9 (x + y)^2 - 15 x (x + y)^2 \\
    &\left.-15 y (x + y)^{2} + 6 {(x + y)}^{3}\right) \\
    &=x^5 + 5 x^4 y + 10 x^3 y^2 + 10 x^2 y^3 + 5 x y^4 + y^5 + 1
\end{align*}
Above example could be verified using using the commands defined in Mathematica package at~\cite{github_source_files}
\begin{itemize}
    \item \texttt{BinomialTheoremAndDiscreteConvolutionTest[0, x + y]}
    \item \texttt{BinomialTheoremAndDiscreteConvolutionTest[1, x + y]}
    \item \texttt{Expand[BinomialTheoremAndDiscreteConvolutionTest[1, x + y]]}
    \item \texttt{BinomialTheoremAndDiscreteConvolutionTest[2, x + y]}
    \item \texttt{Expand[BinomialTheoremAndDiscreteConvolutionTest[2, x + y]]}
\end{itemize}
\begin{cor}
    \label{cor_bin_exp_and_macaulay_conv_strict}
    (Generalization of Theorem~\ref{thm_odd_power_by_macaulays_convolution_strict} for Binomials.)
    For every $m\in\mathbb{N}, \; x,y\in\mathbb{R}$ and $n > 0$
    \begin{equation*}
        \sum_{r=0}^{m} \coeffA{m}{r} \convPower{n}{r}{x+y}
        =
        -1 + \sum_{r=0}^{2m+1} \binom{2m+1}{r} x^{2m+1-r} y^r
    \end{equation*}
\end{cor}
For example, given $m=0,1$ the Corollary~\ref{cor_bin_exp_and_macaulay_conv_strict} gives
\begin{equation*}
    \begin{split}
        \sum_{r=0}^{0} \coeffA{0}{r} \convPower{n}{r}{x+y}
        &= x + y - 1 \\
        \sum_{r=0}^{1} \coeffA{1}{r} \convPower{n}{r}{x+y}
        &= -1 + x + y - (-1 + x + y) (x + y) (-1 - 3 x - 3 y + 2 (x + y)) \\
        &= x^3 + 3 x^2 y + 3 x y^2 + y^3 - 1
    \end{split}
\end{equation*}
Above example could be verified using using the commands defined in Mathematica package at~\cite{github_source_files}
\begin{itemize}
    \item \texttt{BinomialTheoremAndDiscreteConvolutionStrictTest[0, x + y]}
    \item \texttt{BinomialTheoremAndDiscreteConvolutionStrictTest[1, x + y]}
    \item \texttt{Expand[BinomialTheoremAndDiscreteConvolutionStrictTest[1, x + y]]}
\end{itemize}
From the other prospective, the following binomial holds.
For every $n \geq 0$
\begin{equation}
    \label{eq:parametric-identity}
    \begin{split}
        (x-2a)^{2m+1} + 1 &= \sum_{r=0}^{m} \coeffA{m}{r} ((t-k)^r \ast (t-k)^r)[x] \\
                          &= \sum_{r=0}^{m} \coeffA{m}{r} \sum_{k=a}^{x-a} (k-a)^r (x-k-a)^r
    \end{split}
\end{equation}
Similarly, the following binomial holds.
For every $n > 0$
\begin{equation}
    \label{eq:parametric-identity-strict}
    \begin{split}
        (x-2a)^{2m+1} - 1 &= \sum_{r=0}^{m} \coeffA{m}{r} ((t-k)^r \ast (t-k)^r)[x] \\
                          &= \sum_{r=0}^{m} \coeffA{m}{r} \sum_{k=a+1}^{x-a-1} (k-a)^r (x-k-a)^r
    \end{split}
\end{equation}
To validate equations~\eqref{eq:parametric-identity} and~\eqref{eq:parametric-identity-strict}
use the following commands
\begin{itemize}
    \item \texttt{ConvolutionOfBinomial[10, 2, 1]} verifies an equation~\eqref{eq:parametric-identity}.
    \item \texttt{ConvolutionOfBinomial1[10, 2, 1]} verifies an equation~\eqref{eq:parametric-identity-strict}.
\end{itemize}

\subsection{Generalization for Multinomials} \label{subsec:generalization-for-multinomials}
In this subsection we generalize
Theorems~\eqref{thm_odd_power_by_macaulays_convolution} and~\eqref{thm_odd_power_by_macaulays_convolution_strict}
for multinomial cases.
\begin{cor}
    \label{cor_mult_exp_and_macaulay_conv}
    (Generalization of Theorem~\ref{thm_odd_power_by_macaulays_convolution} for Multinomials.)
    For every $x_1, x_2, \ldots, x_t\in\mathbb{R}, \; m\in\mathbb{N}, \; n \geq 1$
    \[
        \sum_{r=0}^{m} \coeffA{m}{r} \convPower{n}{r}{\multifoldSum{t}} =
        1 + \sum_{\multifoldSum[k]{t}=2m+1} \binom{2m+1}{k_1, k_2,\ldots, k_t} \prod_{\ell=1}^{t} x_\ell^{k_\ell}
    \]
\end{cor}
For instance, given $m=1$ the Corollary~\ref{cor_mult_exp_and_macaulay_conv} gives
\begin{equation*}
    \begin{split}
        &\sum_{r=0}^{1} \coeffA{1}{r} \convPower{n}{r}{x+y+z} \\
        &=1 + x + y + z - (x + y + z) (1 + x + y + z) (1 - 3 x - 3 y - 3 z + 2 (x + y + z)) \\
        &=1 + x^3 + 3 x^2 y + 3 x y^2 + y^3 + 3 x^2 z + 6 x y z + 3 y^2 z + 3 x z^2 + 3 y z^2 + z^3.
    \end{split}
\end{equation*}
Above example could be verified using using the commands defined in Mathematica package at~\cite{github_source_files}
\begin{itemize}
    \item \texttt{BinomialTheoremAndDiscreteConvolutionTest[1, x + y + z]}
    \item \texttt{Expand[BinomialTheoremAndDiscreteConvolutionTest[1, x + y + z]]}
\end{itemize}
\begin{cor}
    \label{cor_mult_exp_and_macaulay_conv_strict}
    (Generalization of Theorem~\ref{thm_odd_power_by_macaulays_convolution_strict} for Multinomials.)
    For each $\multifoldSum{t} \geq 1, \; x_1,x_2,\ldots,x_t\in\mathbb{R}, \; m\in\mathbb{N}, \; n\geq1$
    \[
        \sum_{r=0}^{m} \coeffA{m}{r} \convPower{n}{r}{\multifoldSum{t}} =
        -1 + \sum_{\multifoldSum[k]{t}=2m+1} \binom{2m+1}{k_1, k_2,\ldots, k_t} \prod_{\ell=1}^{t} x_\ell^{k_\ell}
    \]
\end{cor}
For example, given $m=1$ the Corollary~\ref{cor_mult_exp_and_macaulay_conv_strict} gives
\begin{equation*}
    \begin{split}
        &\sum_{r=0}^{1} \coeffA{1}{r} \convPower{n}{r}{x+y+z} \\
        &=-1 + x + y + z - (-1 + x + y + z) (x + y + z) (-1 - 3 x - 3 y - 3 z + 2 (x + y + z)) \\
        &=-1 + x^3 + 3 x^2 y + 3 x y^2 + y^3 + 3 x^2 z + 6 x y z + 3 y^2 z + 3 x z^2 + 3 y z^2 + z^3.
    \end{split}
\end{equation*}
Above example could be verified using using the commands defined in Mathematica package at~\cite{github_source_files}
\begin{itemize}
    \item \texttt{BinomialTheoremAndDiscreteConvolutionStrictTest[1, x + y + z]}
    \item \texttt{Expand[BinomialTheoremAndDiscreteConvolutionStrictTest[1, x + y + z]]}
\end{itemize}

    \section{Derivation of the coefficient \texorpdfstring{$\coeffA{m}{r}$}{A[m,r]}}
    \label{sec:derivation-of-coefficients-a}
    By Lemma~\ref{lemma_polynomial_p_and_odd_power} for every $m\in\mathbb{N}, \; n\in\mathbb{R}$
\begin{equation}
    \label{eq:current_a_def}
    n^{2m+1} = \sum_{r=0}^{m} \coeffA{m}{r} \sum_{k=0}^{n-1} k^r (n-k)^r
\end{equation}
The $\coeffA{m}{r}$ might be evaluated using binomial expansion of $\sum_{k=0}^{n-1} k^r (n-k)^r$
\begin{equation*}
    \sum_{k=0}^{n-1} k^r (n-k)^r
    =\sum_{k=0}^{n-1} k^r \sum_{j=0}^{r} (-1)^j \binom{r}{j} n^{r-j} k^{j}
    =\sum_{j=0}^{r} (-1)^j \binom{r}{j} n^{r-j} \sum_{k=0}^{n-1} k^{r+j}
\end{equation*}
Using Faulhaber's formula $\sum_{k=1}^{n} k^{p} = \frac{1}{p+1}\sum_{j=0}^{p} \binom{p+1}{j}
\bernoulli{j} n^{p+1-j}$ we get
\begin{equation}
    \label{eq:proof1}
    \begin{split}
        \sum_{k=0}^{n-1} k^r (n-k)^r
        &=\sum_{j=0}^{r} \binom{r}{j} n^{r-j} \frac{(-1)^j}{r+j+1}
        \left[\sum_{s} \binom{r+j+1}{s} \bernoulli{s} n^{r+j+1-s} - \bernoulli{r+j+1} \right] \\
        &=\sum_{j,s} \binom{r}{j} \frac{(-1)^j}{r+j+1} \binom{r+j+1}{s} \bernoulli{s} n^{2r+1-s}
        -\sum_{j} \binom{r}{j} \frac{(-1)^j}{r+j+1} \bernoulli{r+j+1} n^{r-j} \\
        &=\sum_{s} \underbrace{\sum_{j} \binom{r}{j} \frac{(-1)^j}{r+j+1} \binom{r+j+1}{s}}_{S(r)}
        \bernoulli{s} n^{2r+1-s} \\
        &-\sum_{j} \binom{r}{j} \frac{(-1)^j}{r+j+1} \bernoulli{r+j+1} n^{r-j}
    \end{split}
\end{equation}
where $\bernoulli{s}$ are Bernoulli numbers and $\bernoulli{1}=\frac{1}{2}$.
Now, we notice that
\begin{equation*}
    \sum_{j} \binom{r}{j} \frac{(-1)^j}{r+j+1} \binom{r+j+1}{s}
    =\begin{cases}
         \frac{1}{(2r+1) \binom{2r}r}, & \text{if } s=0;\\
         \frac{(-1)^r}{s} \binom{r}{2r-s+1}, & \text{if } s>0.
    \end{cases}
\end{equation*}
In particular, the last sum is zero for $0<s\leq r$.
Therefore, expression~\eqref{eq:proof1} takes the form
\begin{equation*}
    \begin{split}
        \sum_{k=0}^{n-1} k^r (n-k)^r
        &=\frac{1}{(2r+1) \binom{2r}{r}} n^{2r+1}
        +\underbrace{\sum_{s \geq 1} \frac{(-1)^r}{s} \binom{r}{2r-s+1} \bernoulli{s} n^{2r+1-s}}_{(\star)} \\
        &-\underbrace{\sum_{j} \binom{r}{j} \frac{(-1)^j}{r+j+1} \bernoulli{r+j+1} n^{r-j}}_{(\diamond)}
    \end{split}
\end{equation*}
Hence, by introducing $\ell=2r+1-s$ into $(\star)$ and $\ell=r-j$ into $(\diamond)$, we get
\begin{equation*}
    \begin{split}
        \sum_{k=0}^{n-1} k^r (n-k)^r
        &=\frac{1}{(2r+1) \binom{2r}{r}} n^{2r+1}
        +\sum_{\ell} \frac{(-1)^r}{2r+1-\ell} \binom{r}{\ell} \bernoulli{2r+1-\ell} n^{\ell} \\
        &-\sum_{\ell} \binom{r}{\ell} \frac{(-1)^{j-\ell}}{2r+1-\ell} \bernoulli{2r+1-\ell} n^{\ell}
    \end{split}
\end{equation*}
\begin{equation*}
    \begin{split}
        \sum_{k=0}^{n-1} k^r (n-k)^r
        &=\frac{1}{(2r+1) \binom{2r}{r}} n^{2r+1}
        +(-1)^{r} \sum_{\ell} \frac{1}{2r+1-\ell} \binom{r}{\ell} \bernoulli{2r+1-\ell} n^{\ell} \\
        &-\frac{1}{(-1)^{r}} \sum_{\ell} \binom{r}{\ell} \frac{(-1)^{j-\ell}}{2r+1-\ell} \bernoulli{2r+1-\ell} n^{\ell} \\
        &=\frac{1}{(2r+1) \binom{2r}{r}}n^{2r+1}
        +2 \sum_{\text{odd } \ell}^{r} \frac{(-1)^r}{2r+1-\ell} \binom{r}{\ell} \bernoulli{2r+1-\ell} n^{\ell}
    \end{split}
\end{equation*}
Using the definition~\eqref{eq:current_a_def} of $\coeffA{m}{r}$, we obtain the following identity for polynomials in $n$
\begin{equation}
    \label{eq:proof2}
    \sum_{r=0}^{m} \coeffA{m}{r} \frac{1}{(2r+1) \binom{2r}{r}} n^{2r+1}
    +2 \sum_{r=0}^{m}\sum_{\text{odd } \ell}^{r} \coeffA{m}{r} \frac{(-1)^r}{2r+1-\ell}
    \binom{r}{\ell} \bernoulli{2r+1-\ell} n^{\ell}
    \equiv
    n^{2m+1}
\end{equation}
Taking the coefficient of $n^{2r+1}$ for $r=m$ in~\eqref{eq:proof2} we get $\coeffA{m}{m} = (2m+1) \binom{2m}{m}$.
Since that $\text{odd } \ell \leq r$ in explicit form is $2j + 1 \leq r$, it follows that $j \leq \frac{m-1}{2}$,
where $j$ is an iterator.
Therefore, taking the coefficient of $n^{2j+1}$ for an integer $j$ in the range $\frac{m}{2} \leq j \leq m$,
we get $\coeffA{m}{j} = 0$.
Taking the coefficient of $n^{2d+1}$ for $d$ in the range $m/4 \leq d < m/2$ we get
\begin{equation*}
    \coeffA{m}{d} \frac{1}{(2d+1) \binom{2d}{d}}
    +2 (2m+1) \binom{2m}{m} \binom{m}{2d+1} \frac{(-1)^m}{2m-2d} \bernoulli{2m-2d} = 0,
\end{equation*}
i.e
\begin{equation*}
    \coeffA{m}{d} = (-1)^{m-1} \frac{(2m+1)!}{d!d!m!(m-2d-1)!} \frac{1}{m-d} \bernoulli{2m-2d}
\end{equation*}
Continue similarly we can express $\coeffA{m}{r}$ for each integer $r$ in range $m/2^{s+1}\leq r < m/2^s$
(iterating consecutively $s=1,2,\ldots$) via previously determined values of $\coeffA{m}{d}$ as follows
\begin{equation*}
    \coeffA{m}{r} =
    (2r+1) \binom{2r}{r} \sum_{d=2r+1}^{m} \coeffA{m}{d} \binom{d}{2r+1} \frac{(-1)^{d-1}}{d-r}
    \bernoulli{2d-2r}
\end{equation*}
So that
\begin{align*}
    \coeffA{m}{r} =
    \begin{cases}
    (2r+1)
        \binom{2r}{r} & \text{if } r=m \\
        (2r+1) \binom{2r}{r} \sum_{d=2r+1}^{m} \coeffA{m}{d} \binom{d}{2r+1} \frac{(-1)^{d-1}}{d-r}
        \bernoulli{2d-2r} & \text{if } 0 \leq r<m \\
        0 & \text{if } r<0 \text{ or } r>m
    \end{cases}
\end{align*}
As desired.

    \section{Conclusion}
    \label{sec:conclusion}
    In this manuscript, we introduced the polynomial $\mathbf{P}^{m}_{b}(x)$ and examined its properties.
We established a polynomial identity for odd-powers that demonstrates the connection between Binomial theorem
and discrete convolution of odd-powered polynomials.
This relationship was extended to the multinomial case.
All results were verified using Mathematica programs.

    \section{Acknowledgements}
    \label{sec:acknowledgements}
    I'd like to thank to Dr. Max Alekseyev for sufficient help in the derivation of the real coefficients $\coeffA{m}{r}$.
Also, I'd like to thank to OEIS editors Michel Marcus, Peter Luschny, Jon E. Schoenfield and others
for their useful volunteer work and for useful comments during the work on OEIS sequences related to this manuscript.

    \bibliographystyle{unsrt}
    \bibliography{OnTheBinomialTheoremAndDiscreteConvolutionReferences}

\begin{thebibliography}{10}

\bibitem{WeissteinBernoulli}
Eric~W Weisstein.
\newblock {"Bernoulli Number." From MathWorld -- A Wolfram Web Resource.}
\newblock \url {http://mathworld.wolfram.com/BernoulliNumber.html}.

\bibitem{damelin_discrete_convolution}
Steven B.~Damelin and Willard Miller.
\newblock The mathematics of signal processing.
\newblock {\em The Mathematics of Signal Processing}, page 232, 01 2011.

\bibitem{AbraSteg72}
{Abramowitz, Milton and Stegun, Irene A.}, editor.
\newblock {\em {Handbook of Mathematical Functions with Formulas, Graphs, and
  Mathematical Tables}}.
\newblock U.S. Government Printing Office, Washington, DC, USA, tenth printing
  edition, 1972.

\bibitem{kolosov2018numerator}
Petro Kolosov.
\newblock {Entry A302971 in The On-Line Encyclopedia of Integer Sequences}.
\newblock Published electronically at \url{https://oeis.org/A302971}, 2018.

\bibitem{kolosov2018denominator}
Petro Kolosov.
\newblock {Entry A304042 in The On-Line Encyclopedia of Integer Sequences}.
\newblock Published electronically at \url{https://oeis.org/A304042}, 2018.

\bibitem{kolosov2017third}
Petro Kolosov.
\newblock {Numerical triangle, row sums give third power, Entry A287326 in The
  On-Line Encyclopedia of Integer Sequences}.
\newblock Published electronically at \url{https://oeis.org/A287326}, 2017.

\bibitem{kolosov2018fifth}
Petro Kolosov.
\newblock {Numerical triangle, row sums give fifth power, Entry A300656 in The
  On-Line Encyclopedia of Integer Sequences}.
\newblock Published electronically at \url{https://oeis.org/A300656}, 2018.

\bibitem{iverson_apl}
Kenneth~E. Iverson.
\newblock A programming language.
\newblock In {\em Proceedings of the May 1-3, 1962, Spring Joint Computer
  Conference}, AIEE-IRE '62 (Spring), pages 345--351, New York, NY, USA, 1962.
  ACM.

\bibitem{knuth_two_notes_on_notation}
Donald~E. Knuth.
\newblock {Two Notes on Notation}.
\newblock {\em Am. Math. Monthly}, 99(5):403--422, 1992.

\bibitem{graham1994concrete}
{Graham, Ronald L. and Knuth, Donald E. and Patashnik, Oren}.
\newblock {\em {Concrete mathematics: A foundation for computer science (second
  edition)}}.
\newblock {Addison-Wesley Publishing Company, Inc.}, 1994.
\newblock \url {https://archive.org/details/concrete-mathematics}.

\bibitem{github_source_files}
Petro Kolosov.
\newblock {On the link between binomial theorem and discrete convolution --
  Source files}, 2022.
\newblock \url
  {https://github.com/kolosovpetro/OnTheBinomialTheoremAndDiscreteConvolution}.

\end{thebibliography}
    \noindent \textbf{Version:} \texttt{1.0.5-tags-v1-0-4.2+tags/v1.0.4.5291b99}

    \clearpage

    \section{Addendum 1: Verification of the results}
    \label{sec:verification-of-the-results-and-examples}
    To fulfill our study we provide an opportunity to verify its results by means of Wolfram Mathematica language.

\subsection{Mathematica commands} \label{subsec:mathematica-commands}
Proceeding to the repository~\cite{github_source_files} reader is able to find there a folder named \texttt{mathematica}
that contains the files
\begin{itemize}
    \item \texttt{OnTheBinomialTheoremAndDiscreteConvolution.m} is a package file with definitions
    \item \texttt{OnTheBinomialTheoremAndDiscreteConvolution.nb} is a notebook file with examples.
\end{itemize}
The following commands may be used to reproduce the results of this manuscript:
\begin{itemize}
    \item \texttt{A[m, r]} returns the real coefficient $\coeffA{m}{r}$ defined by~\eqref{eq:def_coeff_a}.
    \item \texttt{PrintTriangleOfA[rows]} prints the table of coefficients $\coeffA{m}{r}$. \\
    Command \texttt{PrintTriangleOfA[7]} reproduces the table (\ref{tab:table_of_coefficients_a}).
    \item \texttt{PolynomialL[m, n, k]} returns the polynomial $\polynomialL{m}{n}{k}$ defined by~\eqref{eq:def_polynomial_l}.
    \item \texttt{PolynomialP[m, x, b]} returns the polynomial $\polynomialP{m}{b}{x}$ defined by~\eqref{eq:def_polynomial_p}.
    \item \texttt{Expand[PolynomialP[m, x + y, x + y]]} verifies the Lemma~\ref{lemma_polynomial_p_and_odd_power}.
    \item \texttt{PolynomialH[m, t, j]} returns the polynomial $\coeffH{m}{t}{j}$ defined by~\eqref{eq:def_coeff_h}.
    \item \texttt{PolynomialX[m, t, k]} returns the polynomial $\polynomialX{m}{t}{k}$ defined by~\eqref{eq:def_coeff_x}.
    \item \texttt{Expand[BinomialTheoremAndDiscreteConvolutionTest[m, x + y]]} verifies the Corollary~\ref{cor_bin_exp_and_macaulay_conv}.
    \item \texttt{Expand[BinomialTheoremAndDiscreteConvolutionStrictTest[m, x + y]]} verifies the Corollary~\ref{cor_bin_exp_and_macaulay_conv_strict}.
    \item \texttt{DiscreteConvolutionPowerIdentityParametricTest[m, x, a]} verifies an equation~\eqref{eq:parametric-identity}.
    Usage \texttt{Column[Table[DiscreteConvolutionPowerIdentityParametricTest[1, x, 1], {x, 3, 20}], Left]}.
    \item \texttt{DiscreteConvolutionPowerIdentityStrictParametricTest[m, x, a]} verifies an equation~\eqref{eq:parametric-identity-strict}.
    Usage \texttt{Column[Table[DiscreteConvolutionPowerIdentityStrictParametricTest[1, x, 1], {x, 3, 20}], Left]}.
    \item \texttt{Expand[PolynomialIdentityOfP[1, n, b]]} validates an identity
    \[\polynomialP{m}{b}{x} = \sum_{r=0}^{m} \coeffA{m}{r} \sum_{j=0}^{r} (-1)^{j} x^{r-j} \binom{r}{j} \sum_{k=0}^{b-1} k^{r+j}\]
    \item \texttt{PolynomialIdentityInvolvingX[m, x, b]} validates an identity~\eqref{eq:p_all_forms}
    \[\polynomialP{m}{b}{x} = \sum_{r=0}^{m} (-1)^{m-r} \polynomialX{m}{r}{b} \cdot x^r\]
    \item \texttt{PolynomialIdentityInvolvingH[m, n, b]} validates an identity~\eqref{eq:p_all_forms}.
    \[\polynomialP{m}{b}{x} =\sum_{r=0}^{m} \sum_{\ell=1}^{2m-r+1} (-1)^{2m-r} \coeffH{m}{r}{\ell} \cdot b^\ell \cdot x^r\]
\end{itemize}

\subsection{Examples} \label{subsec:examples}
For example, given $m=1$ we have the following values of $\polynomialL{1}{x}{k}$
\begin{table}[H]
    \setlength\extrarowheight{-6pt}
    \begin{tabular}{c|cccccccc}
        $x/k$ & 0 & 1  & 2  & 3  & 4  & 5  & 6  & 7 \\
        \hline
        0     & 1 &    &    &    &    &    &    &   \\
        1     & 1 & 1  &    &    &    &    &    &   \\
        2     & 1 & 7  & 1  &    &    &    &    &   \\
        3     & 1 & 13 & 13 & 1  &    &    &    &   \\
        4     & 1 & 19 & 25 & 19 & 1  &    &    &   \\
        5     & 1 & 25 & 37 & 37 & 25 & 1  &    &   \\
        6     & 1 & 31 & 49 & 55 & 49 & 31 & 1  &   \\
        7     & 1 & 37 & 61 & 73 & 73 & 61 & 37 & 1
    \end{tabular}~\caption{Values of $\polynomialL{1}{x}{k}$.
    See OEIS entry: \href{https://oeis.org/A300656}{\texttt{A300656}}, \cite{kolosov2017third}.}
    \label{tab:tab_3}
\end{table}
Table~\ref{tab:tab_3} can be reproduced using Mathematica command
\begin{center}
    \texttt{PrintTriangleOfPolynomialL[1, 7]}
\end{center}
defined in the~\cite{github_source_files}.
From Table~\ref{tab:tab_3} it is seen that
\begin{equation*}
    \begin{split}
        \polynomialP{1}{0}{0} &= 0 = 0^3 \\
        \polynomialP{1}{1}{1} &= 1 = 1^3 \\
        \polynomialP{1}{2}{2} &= 1+7 = 2^3 \\
        \polynomialP{1}{3}{3} &= 1+13+13 = 3^3 \\
        \polynomialP{1}{4}{4} &= 1+19+25+19 = 4^3 \\
        \polynomialP{1}{5}{5} &= 1+25+37+37+25 = 5^3
    \end{split}
\end{equation*}
Another case, given $m=2$ we have the following values of $\polynomialL{2}{x}{k}$
\begin{table}[H]
    \setlength\extrarowheight{-6pt}
    \begin{tabular}{c|cccccccc}
        $x/k$ & 0 & 1    & 2    & 3    & 4    & 5    & 6    & 7 \\
        \hline
        0     & 1 &      &      &      &      &      &      &   \\
        1     & 1 & 1    &      &      &      &      &      &   \\
        2     & 1 & 31   & 1    &      &      &      &      &   \\
        3     & 1 & 121  & 121  & 1    &      &      &      &   \\
        4     & 1 & 271  & 481  & 271  & 1    &      &      &   \\
        5     & 1 & 481  & 1081 & 1081 & 481  & 1    &      &   \\
        6     & 1 & 751  & 1921 & 2431 & 1921 & 751  & 1    &   \\
        7     & 1 & 1081 & 3001 & 4321 & 4321 & 3001 & 1081 & 1
    \end{tabular}
    \caption{Values of $\polynomialL{2}{x}{k}$.
    See the OEIS entry \href{https://oeis.org/A300656}{\texttt{A300656}}, ~\cite{kolosov2018fifth}.}
    \label{tab:tab_4}
\end{table}
Table~\ref{tab:tab_4} can be reproduced using Mathematica command
\begin{center}
    \texttt{PrintTriangleOfPolynomialL[2, 7]}
\end{center}
defined in the~\cite{github_source_files}.
Again, an odd-power identity~\ref{lemma_polynomial_p_and_odd_power} holds
\begin{equation*}
    \begin{split}
        \polynomialP{2}{0}{0} &= 0 = 0^5 \\
        \polynomialP{2}{1}{1} &= 1 = 1^5 \\
        \polynomialP{2}{2}{2} &= 1+31 = 2^5 \\
        \polynomialP{2}{3}{3} &= 1+121+121 = 3^5 \\
        \polynomialP{2}{4}{4} &= 1+271+481+271 = 4^5 \\
        \polynomialP{2}{5}{5} &= 1+481+1081+1081+481 = 5^5
    \end{split}
\end{equation*}

\end{document}